\documentclass[graybox]{svmult}

\usepackage{mathptmx}       
\usepackage{helvet}         
\usepackage{courier}        

\usepackage{amsmath, amssymb, mathtools}
\usepackage{cite}
\usepackage{color,soul}  
\usepackage{flafter}  
\usepackage{graphicx}
\usepackage[mathscr]{eucal}
\usepackage[small]{caption}
\usepackage{subfigure}
\usepackage{bm} 

\newif\ifupdatetikz
\updatetikztrue 
\ifupdatetikz
	\usepackage{tikz}
	\usepackage{pgfplots}
	\pgfplotsset{
			compat=newest,
			tick label style={font=\scriptsize},
			label style={font=\scriptsize},
			title style={font=\scriptsize\bfseries, yshift=-0.2cm},
			legend style={font=\tiny}
		}
	\usepgfplotslibrary{external}
	\usetikzlibrary{shapes,arrows,positioning}
	\tikzset{cross/.style={cross out, draw=black, minimum size=2*(#1-\pgflinewidth), inner sep=0pt, outer sep=0pt},cross/.default={2pt}}
	\tikzstyle{axis} = [-latex',line width=1.25]
	\newcommand{\mygreen}{green!60!black}
\else
	\usepackage{tikzexternal}
	\tikzsetexternalprefix{image/tikz/}
	
\fi
\newcommand{\figname}[1]{\tikzsetnextfilename{#1}}

\newcommand{\datafile}[1]{image/tikz/data/#1.tsv}

\newcommand{\drawDot}[3]
{
	\draw[color=#2,fill=#2] #1 circle (.4ex) #3; 
}

\let\originalleft\left
\let\originalright\right
\renewcommand{\left}{\mathopen{}\mathclose\bgroup\originalleft}
\renewcommand{\right}{\aftergroup\egroup\originalright}

\newcommand{\abs}[1]{\left\vert#1\right\vert}

\newcommand{\cI}{\mathcal{I}}

\newcommand{\cM}{\mathcal{M}}

\DeclareMathOperator{\diverg}{div}

\newcommand{\collOp}{\mathcal{Q}}

\newcommand{\collOper}[1]{\collOp\left(#1\right)(\v)}

\newcommand{\Dtime}{\mathrm{D}_t}
\newcommand{\Dxv}{\mathrm{D}_{\xGrid,\vGrid}}
\newcommand{\Dt}{\Delta t}
\newcommand{\dt}{\delta t}
\newcommand{\dx}{\Delta x}
\newcommand{\fepsi}{f^\epsi}
\newcommand{\fSDSys}{\mathbf{f}} 		

\newcommand{\heatflux}{\mathbf{q}}
\newcommand{\hydroLimit}{\epsi \to 0}
\newcommand{\Ma}{\mathit{Ma}} 

\newcommand{\Maxwellian}[2]{\cM_{#1}\left(#2\right)} 
\newcommand{\MaxArtSDSys}[1]{\Maxwellian{\vGrid}{#1}}	
\newcommand{\MaxBGKOneD}{\Maxwellian{v}{\fepsi}} 		
\newcommand{\MaxBGKMultiD}{\Maxwellian{\v}{\fepsi}}		
\newcommand{\MaxGlobMultiD}{\cM^{\rhoinf,\vMacroMultiDinf,\Tinf}_{\v}} 

\newcommand{\R}{\mathbb{R}}			
\newcommand{\rhoinf}{\rho^{\infty}}

\newcommand{\Sdt}{S_{\dt}}

\newcommand{\sett}[3]{\left(#1\right)_{#2=1}^{#3}}

\newcommand{\Tinf}{T^{\infty}}

\renewcommand{\v}{\mathbf{v}}
\newcommand{\vGrid}{\boldsymbol{v}}
\newcommand{\vMacroOneD}{\bar{v}}

\newcommand{\vMacroMultiD}{\mathbf{\vMacroOneD}}
\newcommand{\vMacroMultiDinf}{\vMacroMultiD^{\infty}}

\newcommand{\vpec}{\mathbf{c}}

\newcommand{\x}{\mathbf{x}} 		
\newcommand{\xGrid}{\boldsymbol{x}}
\newcommand{\epsi}{\varepsilon}		

\newcommand{\Rnm}[2]{\R^{#1 \times #2}}

\begin{document}

\title{Projective integration for nonlinear BGK kinetic equations}
\titlerunning{Projective integration for BGK}

\author{Ward Melis, Thomas Rey and Giovanni Samaey}

\institute{
Ward Melis and Giovanni Samaey
\at 
NUMA (Numerical Analysis and Applied Mathematics), Dept. Computer Science, KU Leuven, \\
Celestijnenlaan 200A, 3001 Leuven, Belgium \\
\email{ward.melis@cs.kuleuven.be, giovanni.samaey@cs.kuleuven.be}
\and 
Thomas Rey \at  
Laboratoire Paul Painlev\'e, Universit\'e de Lille,\\
Cit\'e Scientifique, 59655 Villeneuve d'Ascq, France \\ 
\email{thomas.rey@math.univ-lille1.fr}}

\maketitle

\abstract{	
	We present a high-order, fully explicit, asymptotic-preserving projective integration scheme for the nonlinear BGK equation. The method first takes a few small (inner) steps with a simple, explicit method (such as direct forward Euler) to damp out the stiff components of the solution. Then, the time derivative is estimated and used in an (outer) Runge-Kutta method of arbitrary order. Based on the spectrum of the linearized BGK operator, we deduce that, with an appropriate choice of inner step size, the time step restriction on the outer time step as well as the number of inner time steps is independent of the stiffness of the BGK source term. We illustrate the method with numerical results in one and two spatial dimensions.
  \keywords{Projective integration, BGK, asymptotic-preserving, WENO
  \\[5pt]
  }	
	}

\section{Introduction} \label{sec:introduction}
The Boltzmann equation constitutes the cornerstone of the kinetic theory of rarefied gases. In a dimensionless, scalar setting, it describes the evolution of the one-particle mass distribution function $\fepsi(\x,\v,t) \in \R^{+}$ as:
\begin{equation} \label{eq:Boltzmann_equation}
	\partial_t \fepsi + \v \cdot \nabla_{\x} \fepsi = \frac{1}{\epsi}\collOper{\fepsi},
\end{equation}
where $t \ge 0$ represents time, and $(\x,\v) \subset \Rnm{D_x}{D_v}$ are the $D_x$-dimensional particle positions and $D_v$-dimensional particle velocities. In equation \eqref{eq:Boltzmann_equation}, the dimensionless constant $\epsi > 0$ determines the regime of the gas flow, for which we roughly identify the hydrodynamic regime $(\epsi \le 10^{-4})$, the transitional regime $(\epsi \in [10^{-4},10^{-1}])$, and the kinetic regime $(\epsi \ge 10^{-1})$. Furthermore, the left hand side of \eqref{eq:Boltzmann_equation} corresponds to a linear transport operator that comprises the convection of particles in space, whereas the right hand side contains the Boltzmann collision operator that entails velocity changes due to particle collisions. However, due to its high-dimensional and complicated structure, the Boltzmann collision operator is often replaced by simpler collision models that capture most essential features of the former. The most well-known such model is the BGK model \cite{Bhatnagar1954}, which models collisions as a linear relaxation towards thermodynamic equilibrium, and is given by:
\begin{equation} \label{eq:bgk_equation}
	\partial_t \fepsi + \v \cdot \nabla_{\x} \fepsi = \frac{1}{\epsi}(\MaxBGKMultiD - \fepsi),
\end{equation}
in which $\MaxBGKMultiD$ denotes the local Maxwellian distribution, which, for a $D_v$-dimensional velocity space, is given by:
\begin{equation} \label{eq:maxwellian} 
	\MaxBGKMultiD = \frac{\rho}{(2\pi T)^{D_v/2}} \exp{\left(-\frac{|\v-\vMacroMultiD|^2}{2T}\right)} := \cM_\v^{\rho,\vMacroMultiD,T}.
\end{equation}
The Maxwellian distribution contains the velocity moments of the distribution function $\fepsi$, which are calculated as:
\begin{equation} \label{eq:f_moments} 
	\rho = \int_{\R^{D_v}} \fepsi d\v, \qquad
	\vMacroOneD^d = \frac{1}{\rho}\int_{\R^{D_v}} v^d \fepsi d\v, \qquad
	T = \frac{1}{D_v\rho}\int_{\R^{D_v}} \abs{\v-\vMacroMultiD}^2 \fepsi d\v,
\end{equation}
where $\rho \in \R^+$, $\vMacroMultiD = \sett{\vMacroOneD^d}{d}{D_v} \in \R^{D_v}$ and $T \in \R^+$ are the density, macroscopic velocity and temperature, respectively, which all depend on space $\x$ and time $t$.
Then, in the limit $\hydroLimit$, the solution to equation \eqref{eq:bgk_equation} converges towards $\cM_\v^{\rho,\vMacroMultiD,T}$, whose moments in \eqref{eq:f_moments} are solution to the compressible Euler system:
\begin{equation}
  \label{eq:fluid_euler}
	\left\{ \begin{aligned}
	  & \partial_t \rho  + \diverg_\x  (\rho \, \vMacroMultiD ) \,=\, 0, 				  \\
    & \partial_t(\rho \, \vMacroMultiD ) + \diverg_\x  \left(\rho \, \vMacroMultiD \otimes \vMacroMultiD \,+\, \rho \, T  \,{\rm\bf I}\right) \, =\, \bm{0}, 
				  \\
   & \partial_t E + \diverg_\x \left ( \vMacroMultiD \left ( E +\rho \, T\right )  \right ) \,=\, 0,
	\end{aligned} \right.
\end{equation}
in which $E$ is the second moment of $\fepsi$, namely its total energy.

In this paper, we construct a fully explicit, asymptotic-preserving, arbitrary order time integration method for the stiff equation \eqref{eq:bgk_equation}. For a comprehensive review of numerical schemes for collisional kinetic equations such as equation \eqref{eq:Boltzmann_equation}, we refer to \cite{DimarcoPareschi15}. The asymptotic-preserving property \cite{Jin1999} implies that, in the limit when $\epsi$ tends to zero, an $\epsi$-independent time step constraint, of the form $\Dt = O(\dx)$, can be used, in agreement with the classical hyperbolic CFL constraint for the limiting fluid equations \eqref{eq:fluid_euler}. To achieve this, we will use a projective integration method, which was introduced in  \cite{Gear2003projective} and first applied to kinetic equations in \cite{lafitte2012}.

The remainder of this paper is structured as follows.  We describe the projective integration method in more detail in section \ref{sec:projective_integration}, after which we  discuss (in section \ref{sec:spectral_properties}) the spectral properties of the linearized BGK operator, which are needed to ensure stability of the method. Some numerical experiments are done in section \ref{sec:results}.

\section{Projective integration} \label{sec:projective_integration}

\begin{figure}[t]
	\begin{center}
		\figname{PI_sketch} \newcommand{\dtLength}{0.3}
\newcommand{\DtLength}{3.5}
\newcommand{\npts}{3}
\newcommand{\yCo}{{{1.6,0.9,0.79,0.94},{2.2,1.65,1.5,1.46}}}

\begin{tikzpicture}[node distance = 2cm, auto]
	\draw[axis] (-0.2,0) -- (9,0) node[align=center,right] {time};
	\draw[axis] (0,-0.1) -- (0,2.5);
	
	\draw [thick] (0.3,0.3) to [out=35,in=165] (8.5,0.7);
	
	\def\ta {0.5}
	\pgfmathparse{\yCo[0][0]} \pgfmathsetmacro\ya{\pgfmathresult}
	\draw[dotted] (\ta,0.1) -- (\ta,\ya); \draw (\ta,0.1) -- (\ta,-0.1) node[below] {$t^{n-1}$};
	\pgfmathparse{\yCo[0][2]} \pgfmathsetmacro\ya{\pgfmathresult}
	\pgfmathparse{\yCo[1][0]} \pgfmathsetmacro\yb{\pgfmathresult}
	\draw[dashed] (\ta+2*\dtLength,\ya) -- (\ta+\DtLength,\yb);
	\foreach \iy in {0, ..., \npts}
	{
		\pgfmathparse{\yCo[0][\iy]}
		\def\ys {\pgfmathresult}
		\drawDot{(\ta+\iy*\dtLength,\ys)} {black} {};
	}
	
	\def\tb {\ta+\DtLength}
	\pgfmathparse{\yCo[1][0]} \pgfmathsetmacro\yb{\pgfmathresult}
	\draw[dotted] (\tb,0.1) -- (\tb,\yb); \draw (\tb,0.1) -- (\tb,-0.1) node[below] {$t^{n}$};
	\pgfmathparse{\yCo[1][2]} \pgfmathsetmacro\yb{\pgfmathresult}
	\def\yc {1.3}
	\draw[dashed] (\tb+2*\dtLength,\yb) -- (\tb+\DtLength,\yc);
	\foreach \iy in {0, ..., \npts}
	{
		\pgfmathparse{\yCo[1][\iy]}
		\def\ys {\pgfmathresult}
		\drawDot{(\tb+\iy*\dtLength,\ys)} {black} {};
	}
	
	\def\tc {\tb+\DtLength}
	\draw[dotted] (\tc,0.1) -- (\tc,\yc); \draw (\tc,0.1) -- (\tc,-0.1) node[below] {$t^{n+1}$};
	\drawDot{(7.5,\yc)} {black} {};

	\path
	([shift={(-5\pgflinewidth,-5\pgflinewidth)}]current bounding box.south west)
	([shift={( 2\pgflinewidth, 2\pgflinewidth)}]current bounding box.north east);
\end{tikzpicture}
	\end{center}
  	\vspace*{-0.5cm}\caption{\label{fig:proj_int} Sketch of projective integration. At each time, an explicit method is applied over a number of small time steps (black dots) so as to stably integrate the fast modes. As soon as these modes are sufficiently damped the solution is extrapolated using a much larger time step (dashed lines). }
\end{figure}
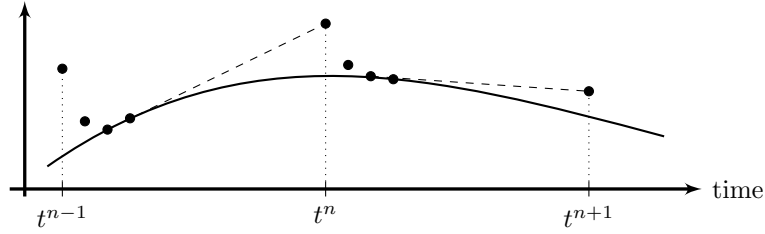

Projective integration \cite{Gear2003projective,lafitte2012} combines a few small time steps with a naive (\emph{inner}) timestepping method (here, a direct forward Euler discretization) with a much larger (\emph{projective, outer}) time step. The idea is sketched in figure \ref{fig:proj_int}.

\vspace*{0.2cm}\noindent\textbf{Inner integrators.} We discretize equation \eqref{eq:bgk_equation} on a uniform, constant in time, periodic spatial mesh with spacing $\dx$, consisting of $I$ mesh points $x_i=i\dx$, ${1 \le i \le I}$, with $I\dx=1$, and a uniform time mesh with time step $\dt$ and discrete time instants $t^k=k\dt$. Furthermore, we discretize velocity space by choosing $J$ discrete components denoted by $\v_j$. The numerical solution on this mesh is denoted by $f_{i,j}^k$, where we have dropped the superscript $\epsi$ on discretized quantities. We then obtain a semidiscrete system of ODEs of the form:
\begin{equation}\label{eq:semidiscrete} 
	\dot{\fSDSys} = \Dtime(\fSDSys),  \qquad
	\Dtime(\fSDSys) = -\Dxv(\fSDSys) + \frac{1}{\epsi}\left(\MaxArtSDSys{\fSDSys} - \fSDSys\right),
\end{equation}
where $\Dxv(\cdot)$ represents a suitable discretization of the convective derivative $\v \cdot \nabla_\x$ (for instance, using upwind differences), and $\fSDSys$ is a vector of size $I \cdot J$.

As inner integrator, we choose the (explicit) forward Euler method with time step $\dt$, for which we will, later on, use the shorthand notation:
\begin{equation} \label{eq:fe_scheme} 
	\fSDSys^{k+1} = \Sdt(\fSDSys^{k}) = \fSDSys^k + \dt\Dtime(\fSDSys^k), \qquad k = 0, 1, \ldots.
\end{equation}

\noindent\textbf{Outer integrators.} In system \eqref{eq:semidiscrete}, the small parameter $\epsi$ leads to the classical time step restriction of the form $\dt = O(\epsi)$ for the inner integrator. However, as $\epsi$ goes to $0$, we obtain the limiting system \eqref{eq:fluid_euler} for which a standard finite volume/forward Euler method only needs to satisfy a stability restriction of the form $\Dt \le C\dx$, with $C$ a constant that depends on the specific choice of the scheme.

In \cite{lafitte2012}, it was proposed to use a projective integration method to accelerate such a brute-force integration; the idea, originating from \cite{Gear2003projective}, is the following. Starting from a computed numerical solution $\fSDSys^n$ at time $t^n=n\Dt$, one first takes $K+1$ \emph{inner} steps of size $\dt$ using~\eqref{eq:fe_scheme}, denoted as 
$\fSDSys^{n,k+1}$,
in which the superscripts $(n,k)$ denote the numerical solution at time ${t^{n,k}=n\Dt +k\dt}$. The aim is to obtain a discrete derivative to be used in the \emph{outer} step to compute $\fSDSys^{n+1} = \fSDSys^{n+1,0}$ via extrapolation in time:
\begin{equation} \label{eq:pfe_scheme}
	\fSDSys^{n+1} = \fSDSys^{n,K+1} + (\Dt - (K + 1)\dt)\frac{\fSDSys^{n,K+1} - \fSDSys^{n,K}}{\dt}.
\end{equation}

Higher-order projective Runge-Kutta (PRK) methods can be constructed by replacing each time derivative evaluation $\mathbf{k}_s$ in a classical Runge-Kutta method by $K+1$ steps of an inner integrator as follows:
\begin{align}
	s = 1 :\;\; & 
	\begin{dcases} 
		\fSDSys^{n,k+1} &= \fSDSys^{n,k} + \dt\Dtime(\fSDSys^{n,k}), \qquad 0 \le k \le K \\ 
		\mathbf{k}_1 &= \dfrac{\fSDSys^{n,K+1} - \fSDSys^{n,K}}{\dt}
	\end{dcases} \label{eq:PRK_stage_1} \\
	2 \le s \le S :\;\; & 
	\begin{dcases} 
		\fSDSys^{n+c_s,0}_s &= \fSDSys^{n,K+1} + (c_s\Dt-(K+1)\dt) \sum_{l=1}^{s-1}\dfrac{a_{s,l}}{c_s} \mathbf{k}_l, \\
		\fSDSys^{n+c_s,k+1}_s &= \fSDSys^{n+c_s,k}_s + \dt\Dtime(\fSDSys^{n+c_s,k}_s), \qquad 0 \le k \le K \\
		\mathbf{k}_s &= \dfrac{\fSDSys^{n+c_s,K+1}_s - \fSDSys^{n+c_s,K}_s}{\dt}
	\end{dcases} \label{eq:PRK_stage_s} \\
	& \fSDSys^{n+1} = \fSDSys^{n,K+1} + (\Dt-(K+1)\dt)\sum_{s=1}^{S}b_s \mathbf{k}_s.
\end{align}
To ensure consistency, the Runge-Kutta matrix $\mathbf{a}=(a_{s,i})_{s,i=1}^S$, weights ${\mathbf{b}=(b_s)_{s=1}^S}$, and nodes $\mathbf{c}=(c_s)_{s=1}^S$ satisfy the conditions $0\le b_s \le 1$ and $0 \le c_s \le 1,$ as well as:
\begin{equation} \label{eq:RK_conditions}
	\sum_{s=1}^Sb_s=1, \qquad \sum_{i=1}^{S-1} a_{s,i} =c_s, \quad 1 \le s \le S. 
\end{equation}

\section{Spectral properties} \label{sec:spectral_properties}

To choose the method parameters (the size of the small and large time steps $\delta t$ and $\Delta t$, as well as the number $K$ of small steps), one needs to analyze the spectrum of the collision operator.  In \cite{Melis2016}, this was done in the hyperbolic scaling for a  system with a linear Maxwellian that serves as a relaxation of a nonlinear hyperbolic conservation law.

By linearizing the Maxwellian \eqref{eq:maxwellian} around the global Maxwellian distribution $\MaxGlobMultiD = \cM_\v^{1,0,1}$, it is shown in \cite[p.206]{Cercignani1988} that the resulting linearized equilibrium can be written as:
\begin{equation} \label{eq:linearized_bgk_maxwellian}
	\cM_\text{lin}(\fepsi)(\x,\v,t) = \sum_{k=0}^{D_v+1} \Psi_k(\v) (\Psi_k, \fepsi)(\x,t),
\end{equation}
in which the scalar product is defined by:
\begin{equation} \label{eq:scalar_product_Hilbert}
	(g,h) = \int_{\R^{D_v}} g(\v)\overline{h(\v)} \frac{1}{(2\pi)^{D_v/2}}\exp\left(\frac{-\abs{\v}^2}{2}\right) d\v.
\end{equation}
Furthermore, the orthonormal set of basis functions $\Psi_k(\v)$ in \eqref{eq:linearized_bgk_maxwellian} are obtained from a straightforward application of the Gram-Schmidt process to the $D_v+1$ collision invariants $(1,\v,\abs{\v}^2)$, yielding:
\begin{equation} \label{eq:Psi_normalized}
	\big(\Psi_0(\v), \ldots, \Psi_{D_v+1}(\v)\big) = \left(1, v^1, ..., v^{D_v}, \frac{\abs{\v}^2 - D_v}{2^{D_v/2}}\right).
\end{equation}
Using the linearized Maxwellian \eqref{eq:linearized_bgk_maxwellian}, the linearized version of the full BGK equation \eqref{eq:bgk_equation} reads:
\begin{equation} \label{eq:linearized_bgk_equation_framework}
	\partial_t \fepsi + \v \cdot \nabla_{\x} \fepsi = -\frac{1}{\epsi}(\cI - \Pi_\text{BGK})\fepsi,
\end{equation}
where $\cI$ denotes the identity operator and $\Pi_\text{BGK}$ is the following rank-$(D_v+2)$ projection operator:
\begin{equation} \label{eq:projection_operator_bgk}
	\Pi_\text{BGK} \fepsi = \sum_{k=0}^{D_v+1} \Psi_k(\v)(\Psi_k, \fepsi).
\end{equation}

This shows that the structure of the linearized Maxwellian \eqref{eq:linearized_bgk_maxwellian} and the linearized BGK projection operator \eqref{eq:projection_operator_bgk} are almost identical to those in \cite{Melis2016}. We can actually view these linear kinetic models as a special simplified case of the linearized BGK equation.
Therefore, it is expected that the construction of stable, asymptotic-preserving projective integration methods for the full BGK equation \eqref{eq:bgk_equation} is practically identical to that in \cite{Melis2016}. In particular, the conclusion is that, when choosing $\delta t=\epsi$, one is able to choose $\Delta t=O(\Delta x)$ and $K$ independent of $\epsi$, resulting in a scheme with computational cost independent of $\epsi$.

\section{Numerical experiments} \label{sec:results}

\noindent\textbf{BGK in 1D.} As a first experiment, we focus on the nonlinear BGK equation \eqref{eq:bgk_equation} in 1D. We consider a Sod-like test case for $x \in [0,1]$ consisting of an initial centered Riemann problem with the following left and right state values:
\begin{equation} \label{eq:bgk_1d_sod}
	\big(\rho_L, \vMacroOneD_L, T_L\big)  = (1, 0, 1), \qquad\quad 
	\big(\rho_R, \vMacroOneD_R, T_R\big) = (0.125, 0, 0.25).
\end{equation}
The initial distribution $\fepsi(x,v,0)$ is then chosen as the Maxwellian \eqref{eq:maxwellian} corresponding to the above initial macroscopic variables. We impose outflow boundary conditions and perform simulations for $t \in [0,0.15]$. 
As velocity space, we take the interval $[-8,8]$, which we discretize on a uniform grid using $J=80$ velocity nodes. In all simulations, space is discretized using the WENO3 spatial discretization with $\dx = 0.01$. Below, we compare solutions for three gas flow regimes: $\epsi = 10^{-1}$ (kinetic regime), $\epsi = 10^{-2}$ (transitional regime) and $\epsi = 10^{-5}$ (fluid regime).

In the kinetic $(\epsi = 10^{-1})$ and transitional $(\epsi = 10^{-2})$ regimes, we compute the numerical solution using the fourth order Runge-Kutta (RK4) time discretization with time step $\dt = 0.1\dx$. 
In the fluid regime $(\epsi = 10^{-5})$, direct integration schemes such as RK4 become too expensive due to a severe time step restriction, which is required to ensure stability of the method. Exploiting that the spectrum of the linearized BGK equation is close to that of the linear kinetic models used in \cite{Melis2016}, see section \ref{sec:spectral_properties}, we construct a projective integration method to accelerate time integration in the fluid regime. As inner integrator, we select the forward Euler time discretization with $\dt = \epsi$. As outer integrator, we choose the fourth-order projective Runge-Kutta (PRK4) method, using $K=2$ inner steps and an outer step of size $\Dt = 0.4\dx$.

\begin{figure}
	\begin{center}
		\figname{bgk1D_epsi125} 
%
%
\definecolor{mycolor2}{rgb}{0.74902,0.00000,0.74902}%
\newcommand{\figWidth}{4.5cm} 
\newcommand{\figHeight}{3cm} 
\newcommand{\figSpacingRight}{1.5cm} 
\newcommand{\figSpacingTop}{1.2cm} 
\begin{tikzpicture}

\begin{axis}[%
width=\figWidth,
height=\figHeight,
at={(0,\figHeight+\figSpacingTop)},
scale only axis,
xmin=0,
xmax=1,
xlabel={$x$},
ymin=0.05,
ymax=1.075,
ytick = {0, 0.2, ..., 1},
ylabel={$\rho(x,t)$},
axis background/.style={fill=white},
title style={font=\scriptsize\bfseries},
title={Density},
clip mode=individual
]
\addplot [color=blue,solid,mark=*,mark size=0.75pt,mark options={solid},forget plot]
  table[]{\datafile{bgk1D1D_epsi124-1}};
\addplot [color=mycolor2,solid,mark=*,mark size=0.75pt,mark options={solid},forget plot]
  table[]{\datafile{bgk1D1D_epsi124-3}};
\addplot [color=\mygreen,solid,mark=*,mark size=0.75pt,mark options={solid},forget plot]
  table[]{\datafile{bgk1D1D_epsi124-4}};
\addplot [color=red,solid,line width=0.75pt,forget plot]
  table[]{\datafile{bgk1D1D_epsi124-2}};
\end{axis}

\begin{axis}[%
width=\figWidth,
height=\figHeight,
at={(\figWidth+\figSpacingRight,\figHeight+\figSpacingTop)},
scale only axis,
xmin=0,
xmax=1,
xlabel={$x$},
ymin=-0.05,
ymax=0.9,
ytick = {0, 0.2, ..., 1.5},
ylabel={$\vMacroOneD(x,t)$},
axis background/.style={fill=white},
title style={font=\scriptsize\bfseries},
title={Velocity},
clip mode=individual
]
\addplot [color=blue,solid,mark=*,mark size=0.75pt,mark options={solid},forget plot]
  table[]{\datafile{bgk1D1D_epsi124-5}};
\addplot [color=mycolor2,solid,mark=*,mark size=0.75pt,mark options={solid},forget plot]
  table[]{\datafile{bgk1D1D_epsi124-7}};
\addplot [color=\mygreen,solid,mark=*,mark size=0.75pt,mark options={solid},forget plot]
  table[]{\datafile{bgk1D1D_epsi124-8}};
\addplot [color=red,solid,line width=0.75pt,forget plot]
  table[]{\datafile{bgk1D1D_epsi124-6}};
\end{axis}

\begin{axis}[%
width=\figWidth,
height=\figHeight,
at={(0,0)},
scale only axis,
xmin=0,
xmax=1,
xlabel={$x$},
ymin=0.15,
ymax=1.1,
ytick = {0, 0.2, ..., 1},
ylabel={$T(x,t)$},
axis background/.style={fill=white},
title style={font=\scriptsize\bfseries},
title={Temperature},
clip mode=individual
]
\addplot [color=blue,solid,mark=*,mark size=0.75pt,mark options={solid},forget plot]
  table[]{\datafile{bgk1D1D_epsi124-9}};
\addplot [color=mycolor2,solid,mark=*,mark size=0.75pt,mark options={solid},forget plot]
  table[]{\datafile{bgk1D1D_epsi124-11}};
\addplot [color=\mygreen,solid,mark=*,mark size=0.75pt,mark options={solid},forget plot]
  table[]{\datafile{bgk1D1D_epsi124-12}};
\addplot [color=red,solid,line width=0.75pt,forget plot]
  table[]{\datafile{bgk1D1D_epsi124-10}};
\end{axis}

\begin{axis}[%
width=\figWidth,
height=\figHeight,
at={(\figWidth+\figSpacingRight,0)},
scale only axis,
xmin=0,
xmax=1,
xlabel={$x$},
ymin=-0.05,
ymax=0.1,
ytick = {-0.05, 0, ..., 0.1},
yticklabel style={/pgf/number format/fixed},
ylabel={$q(x,t)$},
ylabel style={yshift=-0.2cm},
axis background/.style={fill=white},
title style={font=\scriptsize\bfseries},
title={Heat flux},
clip mode=individual
]
\addplot [color=blue,solid,mark=*,mark size=0.75pt,mark options={solid},forget plot]
  table[]{\datafile{bgk1D1D_epsi124-13}};
\addplot [color=mycolor2,solid,mark=*,mark size=0.75pt,mark options={solid},forget plot]
  table[]{\datafile{bgk1D1D_epsi124-15}};
\addplot [color=\mygreen,solid,mark=*,mark size=0.75pt,mark options={solid},forget plot]
  table[]{\datafile{bgk1D1D_epsi124-16}};
\addplot [color=red,solid,line width=0.75pt,forget plot]
  table[]{\datafile{bgk1D1D_epsi124-14}};
\end{axis}
\end{tikzpicture}%
	\end{center}
	\vspace*{-0.5cm}\caption{\label{fig:bgk_1d} Numerical solution of the BGK equation in 1D at $t = 0.15$ for a Sod-like shock test \eqref{eq:bgk_1d_sod} using WENO3 with $\dx = 0.01$. RK4 is used for $\epsi = 10^{-1}$ (blue dots) and $\epsi = 10^{-2}$ (purple dots). The PRK4 method is used for $\epsi = 10^{-5}$ (green dots). Red line: hydrodynamic limit $(\hydroLimit)$. }
\end{figure}
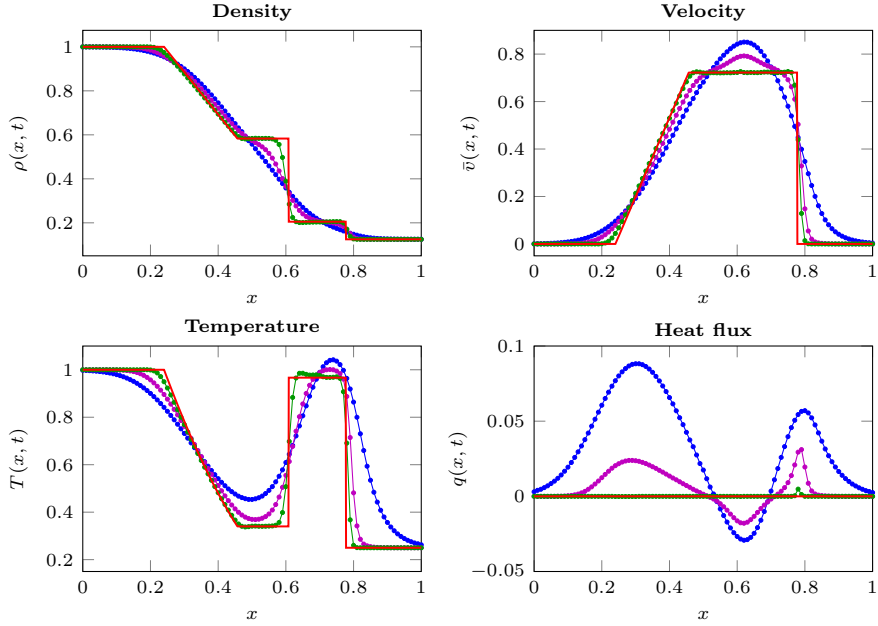

The results are shown in figure \ref{fig:bgk_1d}, where we display the density $\rho$, macroscopic velocity $\vMacroOneD$ and temperature $T$ as given in \eqref{eq:f_moments} at $t = 0.15$. In addition, we plot the heat flux $q$, which, in a general $D_v$-dimensional setting, is a vector $\heatflux = \sett{q^d}{d}{D_v}$ with components given by:
\begin{equation} \label{eq:heat_flux}
	q^d = \frac{1}{2}\int_{\R^{D_v}} \abs{\vpec}^2c^d\fepsi d\v,
\end{equation}
in which $\vpec = \sett{c^d}{d}{D_v} = \v - \vMacroMultiD$ is the peculiar velocity. The different regimes are shown by blue (kinetic), purple (transitional) and green (fluid) dots. The red line in each plot denotes the limiting $(\hydroLimit)$ solution of each macroscopic variable, which all converge to the solution of the compressible Euler equations \eqref{eq:fluid_euler} with ideal gas law $P = \rho T$ and heat flux $q = 0$. From this, we observe that the BGK solution is increasingly dissipative for increasing values of $\epsi$ since the rate with which $\fepsi$ converges to its equilibrium $\MaxBGKOneD$ becomes slower. In contrast, for sufficiently small $\epsi$, relaxation to thermodynamic equilibrium occurs practically instantaneous and the Euler equations \eqref{eq:fluid_euler} yield a valid description. Since this is a hyperbolic system, it allows for the development of sharp discontinuous and shock waves which are clearly seen in the numerical solution.

\vspace*{0.2cm}\noindent\textbf{Shock-bubble interaction in 2D.} Here, we consider the BGK equation in 2D and we investigate the interaction between a moving shock wave and a stationary smooth bubble, which was proposed in \cite{Torrilhon2006}, see also \cite{Cai2010}. This problem consists of a shock wave positioned at $x = -1$ in a spatial domain $\x = (x,y) \in [-2,3] \times [-1,1]$ traveling with Mach number $\Ma = 2$ into an equilibrium flow region. Over the shock wave, the following left $(x \le -1)$ and right $(x > -1)$ state values are imposed \cite{Cai2010}:
\begin{equation} \label{eq:bgk_2d_shock}
	\big(\rho_L, \vMacroOneD^x_L, \vMacroOneD^y_L, T_L\big)  = \left(\frac{16}{7}, \sqrt{\frac{5}{3}}\frac{7}{16}, 0, \frac{133}{64}\right), \qquad\quad
	\big(\rho_R, \vMacroMultiD_R, T_R\big) = \left(1, \bm{0}, 1\right).
\end{equation}
Due to this initial profile, the shock wave will propagate rightwards into the flow region at rest $(x > -1)$. Moreover, in this equilibrium region, a smooth Gaussian density bubble centered at $\x_0 = (0.5,0)$ is placed, given by:
\begin{equation} \label{eq:bgk_2d_bubble}
	\rho(\x,0) = 1 + 1.5\exp\left(-16\abs{\x - \x_0}^2\right).
\end{equation}
Then, the initial distribution $\fepsi(\x,\v,0)$ is chosen as the Maxwellian \eqref{eq:maxwellian} corresponding to the initial macroscopic variables in \eqref{eq:bgk_2d_shock}-\eqref{eq:bgk_2d_bubble}. We impose outflow and periodic boundary conditions along the $x$- and $y$-directions, respectively, and we perform simulations for $t \in [0,0.8]$.
As velocity space, we take the domain $[-10,10]^2$, which we discretize on a uniform grid using $J_x = J_y = 30$. We discretize space using the WENO2 spatial discretization with $I_x = 200$ and $I_y = 25$. Furthermore, we consider a fluid regime by taking $\epsi = 10^{-5}$.

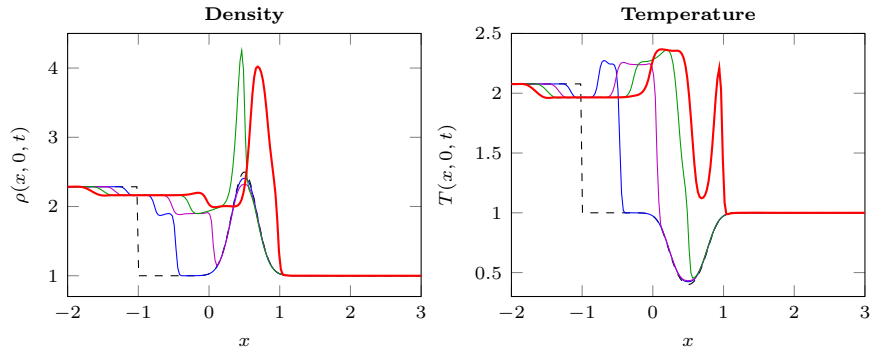
\begin{figure}
	\begin{center}
		\figname{bgk2D_shock_bubble_y0} 
%
%
\definecolor{mycolor2}{rgb}{0.74902,0.00000,0.74902}%
\newcommand{\figWidth}{4.7cm} 
\newcommand{\figHeight}{3.5cm} 
\newcommand{\figSpacingRight}{1.2cm} 
\begin{tikzpicture}

\begin{axis}[%
width=\figWidth,
height=\figHeight,
at={(0,0)},
scale only axis,
xmin=-2.0000,
xmax=3.0000,
xlabel={$x$},
xtick={-2, -1, ..., 3},
ymin=0.7000,
ymax=4.5000,
ylabel={$\rho(x,0,t)$},
axis background/.style={fill=white},
title style={font=\scriptsize\bfseries},
title={Density}
]
\addplot [color=black,dashed,forget plot]
  table[]{\datafile{bgk2D2D_shock_bubble_y0-1}};
\addplot [color=blue,solid,forget plot]
  table[]{\datafile{bgk2D2D_shock_bubble_y0-2}};
\addplot [color=mycolor2,solid,forget plot]
  table[]{\datafile{bgk2D2D_shock_bubble_y0-3}};
\addplot [color=\mygreen,solid,forget plot]
  table[]{\datafile{bgk2D2D_shock_bubble_y0-4}};
\addplot [color=red,solid,line width=0.8pt,forget plot]
  table[]{\datafile{bgk2D2D_shock_bubble_y0-5}};
\end{axis}

\begin{axis}[%
width=\figWidth,
height=\figHeight,
at={(\figWidth+\figSpacingRight,0)},
scale only axis,
xmin=-2.0000,
xmax=3.0000,
xlabel={$x$},
xtick={-2, -1, ..., 3},
ymin=0.3000,
ymax=2.5000,
ylabel={$T(x,0,t)$},
axis background/.style={fill=white},
title style={font=\scriptsize\bfseries},
title={Temperature}
]
\addplot [color=black,dashed,forget plot]
  table[]{\datafile{bgk2D2D_shock_bubble_y0-6}};
\addplot [color=blue,solid,forget plot]
  table[]{\datafile{bgk2D2D_shock_bubble_y0-7}};
\addplot [color=mycolor2,solid,forget plot]
  table[]{\datafile{bgk2D2D_shock_bubble_y0-8}};
\addplot [color=\mygreen,solid,forget plot]
  table[]{\datafile{bgk2D2D_shock_bubble_y0-9}};
\addplot [color=red,solid,line width=0.8pt,forget plot]
  table[]{\datafile{bgk2D2D_shock_bubble_y0-10}};
\end{axis}
\end{tikzpicture}%
	\end{center}
	\vspace*{-0.5cm}\caption{\label{fig:bgk_2d_shock_bubble_y0} Numerical solution of the shock-bubble interaction along $y = 0$ at $t = 0$ (black dashed), $t = 0.2$ (blue), $t = 0.4$ (purple), $t = 0.6$ (green) and $t = 0.8$ (red). }
\end{figure}

We construct a PRK4 method with FE as inner integrator to speed up simulation in time. The inner time step is fixed as $\dt = \epsi$ and we use $K = 2$ inner steps in each outer integrator iteration. The outer time step is chosen as $\Dt = 0.4\dx$. To compare our results with those in \cite{Torrilhon2006}, where the smallest value of $\epsi$ is chosen as $\epsi = 10^{-2}$, we regard the one-dimensional evolution of density and temperature along the axis $y = 0$. For ${t \in \{0, 0.2, 0.4, 0.6, 0.8\}}$, we plot these intersections in figure \ref{fig:bgk_2d_shock_bubble_y0}. We conclude that we obtain the same solution structure at $t = 0.8$ as in \cite{Torrilhon2006}. However, our results are sharper and less dissipative supposedly due to the particular small value of $\epsi$ ($10^{-5}$ versus $10^{-2}$). In contrast to \cite{Cai2010}, we nicely capture the swift changes in the temperature profile for $x \in [0.5,1]$ at $t = 0.8$.

\bibliographystyle{plain} \bibliography{refs}

\end{document}